\documentclass[11pt]{amsart}

\usepackage{hyphenat}

\usepackage[T1]{fontenc}
\usepackage[utf8]{inputenc}
\usepackage{lmodern}
\usepackage{microtype}
\usepackage{amssymb,amsmath,amsthm,mathtools}
\usepackage[colorlinks=true,linkcolor=blue,citecolor=blue,urlcolor=blue]{hyperref}
\usepackage{enumitem}
\usepackage{tikz}
\usepackage{cite}
\usepackage{comment}

\usepackage{minted}
\newmintinline[lean]{lean4}{bgcolor=white}
\newminted[leancode]{lean4}{
  fontsize=\footnotesize,
  breaklines=true,
  bgcolor=white
}

\setlist[itemize]{leftmargin=2.2em}
\numberwithin{equation}{section}

\tikzset{vtx/.style={circle, fill, inner sep=1.5pt}}
\tikzset{openvtx/.style={circle, draw, inner sep=1.5pt}}

\newtheorem{theorem}{Theorem}

\newcommand{\be}{\overline e}
\newcommand{\bT}{\overline T}

\title[Autonomously discovered proof that $R(B_8,B_{10})=37$]{An automated proof that \texorpdfstring{$R(B_8,B_{10})=37$}{R(B8,B10)=37}}

\author{Jeremy Kalfus}
\address{\url{https://jeremykalfus.github.io} Indian Springs School, Indian Springs, AL}
\email{jeremykalfus@gmail.com}

\author{Bernard Lidick\'y}
\address{Department of Mathematics, Iowa State University, Ames, IA, USA.}
\email{lidicky@iastate.edu}
\thanks{Research of this author is supported in part by NSF FRG DMS-2152490, the Simons Foundation TSM-00013439 and Scott Hanna Professorship.}

\date{\today}

\begin{document}

\begin{abstract}
We present a short proof that the book Ramsey number $R(B_8,B_{10})$ equals $37$. The lower bound $R(B_8,B_{10})\ge 37$ is already available in the literature, so it is enough to rule out a $37$-vertex graph containing neither a copy of $B_8$ nor a copy of $B_{10}$ in its complement. 
The problem as well as the proof were found with \emph{AutoMath}, an AI-assisted mathematical discovery workflow developed by the first author.
A Lean formalization of the upper-bound argument is available in the accompanying repository.
\end{abstract}

\maketitle

\section{Introduction}

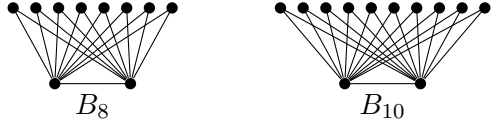
\begin{figure}
\begin{center}
    \begin{tikzpicture}
        \draw (0,0) node[vtx](a){} -- (1,0) node[vtx](b){};
        \foreach \x in {1,...,8}{
        \draw (a) -- (0.5-0.3*4.5+0.3*\x,1) node[vtx]{} -- (b) ;
        }
        \draw (0.5,0) node[below]{$B_8$};
    \end{tikzpicture}
    \hskip 3em
    \begin{tikzpicture}
        \draw (0,0) node[vtx](a){} -- (1,0) node[vtx](b){};
        \foreach \x in {1,...,10}{
        \draw (a) -- (0.5-0.3*5.5+0.3*\x,1) node[vtx]{} -- (b) ;
        }
        \draw (0.5,0) node[below]{$B_{10}$};
    \end{tikzpicture}    
\end{center}
\caption{Book graphs $B_8$ and $B_{10}$.}\label{fig:B8B10}
\end{figure}

The Ramsey number $R(G,H)$ is the smallest integer $n$ such that every graph on $n$ vertices contains $G$ as a subgraph or its complement contains $H$ as a subgraph.  
Results for small graphs $G$ and $H$ are collected in the dynamic survey of Radziszowski~\cite{Rad}.
For $t\ge 1$, let $B_t$ denote the \emph{book graph} consisting of $t$ triangles sharing a common edge. Equivalently, $B_t$ is the graph obtained from an edge $uv$ by adding $t$ further vertices, each adjacent to both $u$ and $v$.
See Figure~\ref{fig:B8B10} for an example.
The study of Ramsey numbers for books was initiated by Rousseau and Sheehan~\cite{RS}, and there have been recent improvements using SAT solvers and related computational methods~\cite{W,LMPV}.
Machine learning, in the form of reinforcement learning, has also been used to find lower bounds for Ramsey numbers involving books~\cite{MSAD}.

Recent work of Lidick\'y, McKinley, Pfender, and Van Overberghe~\cite{LMPV} gives the general lower bound
\[
4n-3\le R(B_{n-2},B_n)\qquad (4\le n\le 21),
\]
so in particular $R(B_8,B_{10})\ge 37$. Before the present argument, the value of $R(B_8,B_{10})$ lay in the one-gap window
\[
37\le R(B_8,B_{10})\le 38,
\]
with the upper bound being a special case of a result of Rousseau and Sheehan~\cite[Theorem 2]{RS}; see also the dynamic survey of Radziszowski~\cite{Rad}. 
In this paper we improve the upper bound to 37, which settles $R(B_8,B_{10})$.

\begin{theorem}
\label{thm:main}
The book Ramsey number $R(B_8,B_{10})$ is equal to $37$.
\end{theorem}

Since the lower bound $R(B_8,B_{10})\ge 37$ is already known, it is enough to prove that every graph on $37$ vertices contains $B_8$ or its complement contains $B_{10}$.

The argument has two parts. The first part uses structural arguments showing that the hypothetical graph on 37 vertices without $B_8$ or $B_{10}$ in its complement is strongly regular. In the second part, the existence of such strongly regular graphs is excluded via a spectral argument. 

In Section~\ref{sec:proof} we describe the proof of Theorem~\ref{thm:main}. 
Section~\ref{sec:lean} is devoted to an explanation of the formalization of Theorem~\ref{thm:main} in Lean.
Finally, Section~\ref{sec:AI} contains a discussion of the software used for this project.

\textbf{AI use.}
The proof of Theorem~\ref{thm:main} was found with \emph{AutoMath}, AI-assisted mathematical-discovery software developed by the first author~\cite{AutoMath}.
The proof was edited for clarity by the authors.
The Lean formalization was produced with GPT-5.5 Pro through Codex.

\section{The proof of Theorem~\ref{thm:main}}\label{sec:proof}
For a graph $G$, denote the number of triangles in $G$ by $t(G)$. 
The complement of $G$ is denoted by $\overline{G}$.

Assume for contradiction that there exists a graph $G=(V,E)$ on $37$ vertices such that $G$ contains no copy of $B_8$ and $\overline G$ contains no copy of $B_{10}$. Write
\[
v=37,\qquad e=|E(G)|,\qquad \be=|E(\overline G)|,\qquad T=t(G),\qquad \bT=t(\overline G).
\]
For two distinct vertices $u,v\in V(G)$, let 
\[
c_G(u,v)=|N_G(u)\cap N_G(v)|
\]
denote the number of their common neighbors.
Notice that $c_{\overline{G}}(u,v)$ is the number of common nonneighbors of $u$ and $v$ in $G$.

Since $G$ has no copy of $B_8$, every edge $uv\in E(G)$ satisfies
\[
c_G(u,v)\le 7.
\]
Analogously, since $\overline G$ contains no $B_{10}$, every edge of $\overline G$ lies in at most $9$ triangles of $\overline G$. 
Hence every nonedge $uv$ of $G$ satisfies
\[
c_{\overline{G}}(u,v)\le 9.
\]
Summing codegrees over all edges gives
\begin{equation}
\label{eq:trianglecaps}
3T=\sum_{uv\in E(G)} c_G(u,v)\le 7e
\qquad
\text{and}
\qquad
3\bT=\sum_{uv\in E(\overline{G})} c_{\overline{G}}(u,v) \le 9\be.
\end{equation}

By combining \eqref{eq:trianglecaps} and $\be=\binom{37}{2}-e=666-e$, we obtain
\[
T+\bT\le \frac73 e+3\be=\frac73 e+3(666-e)=1998-\frac23 e.
\]
Since $2e=\sum_x d(x)$, this implies
\begin{equation}
\label{eq:upperfromcaps}
T+\bT\le 1998-\frac13\sum_{x\in V(G)} d(x).
\end{equation}

Goodman's identity~\cite{Goodman} gives
\begin{equation}
\label{eq:goodman}
T+\bT=\binom{37}{3}-\frac12\sum_{x\in V(G)} d(x)(36-d(x)).
\end{equation}

Combining \eqref{eq:goodman} and \eqref{eq:upperfromcaps} yields
\[
7770-\frac12\sum_x d(x)(36-d(x))\le 1998-\frac13\sum_x d(x).
\]
After rearranging the terms we obtain
\begin{equation}
\label{eq:quadraticineq}
\sum_{x\in V(G)} \bigl(106d(x)-3d(x)^2\bigr)\ge 34632.
\end{equation}

Now we find an upper bound on the left-hand side. It corresponds to sums of the following quadratic polynomial
\[
f(d)=106d-3d^2.
\]
This polynomial is maximized at $d=53/3$. Since $d(x)$ is an integer, the unique maximum is attained at $d=18$, where
\[
f(18)=106\cdot 18-3\cdot 18^2=936.
\]
Therefore
\[
\sum_{x\in V(G)} \bigl(106d(x)-3d(x)^2\bigr)\le 37\cdot 936=34632.
\]
Together with \eqref{eq:quadraticineq}, this forces equality in all of the previous inequalities.
In particular,
\[
d(x)=18\qquad\text{for every }x\in V(G).
\]
Thus $G$ is $18$-regular, and hence
\[
e=\frac{37\cdot 18}{2}=333,
\qquad
\be=666-333=333.
\]

Substituting this degree sequence into Goodman's identity gives
\[
T+\bT=\binom{37}{3}-\frac12\cdot 37\cdot 18\cdot 18
=7770-5994=1776.
\]
Meanwhile the estimate \eqref{eq:trianglecaps} gives
\[
T+\bT\le \frac73\cdot 333+3\cdot 333=777+999=1776.
\]
Hence equality also holds in \eqref{eq:trianglecaps}. Since every summand in
\[
\sum_{uv\in E(G)} c_G(u,v)\le 7e
\]
is an integer at most $7$, equality implies that
\begin{equation}
\label{eq:lambda}
c_G(u,v)=7\qquad\text{for every edge }uv\in E(G).
\end{equation}
Similarly,
\begin{equation}
\label{eq:cn}
c_{\overline{G}}(u,v)=9\qquad\text{for every nonedge }uv\notin E(G).
\end{equation}

Now fix a nonedge $uv$ of $G$, and write $c=c_G(u,v)$. Since $d(u)=d(v)=18$,
\[
|N_G(u)\cup N_G(v)|=18+18-c=36-c.
\]
Among the remaining $35$ vertices, the number adjacent to neither $u$ nor $v$ is therefore
\[
35-|N_G(u)\cup N_G(v)|=35-(36-c)=c-1.
\]
By \eqref{eq:cn}, this number equals $9$, so $c-1=9$ and hence
\begin{equation}
\label{eq:mu}
c_G(u,v)=10\qquad\text{for every nonedge }uv\notin E(G).
\end{equation}
Therefore any hypothetical witness would have to be a strongly regular graph with parameters
\[
(v,k,\lambda,\mu)=(37,18,7,10).
\]

Let $A$ be its adjacency matrix. The standard strongly regular relation is
\[
A^2=kI+\lambda A+\mu(J-I-A).
\]
Substituting $(k,\lambda,\mu)=(18,7,10)$ gives
\begin{equation*}
A^2=18I+7A+10(J-I-A)=-3A+8I+10J.
\end{equation*}
On the orthogonal complement $\mathbf 1^\perp$ of the all-ones vector, the matrix $J$ acts as zero, so every nontrivial eigenvalue $x$ of $A$ satisfies
\[
x^2+3x-8=0.
\]
Thus the two nontrivial eigenvalues are
\[
r,s=\frac{-3\pm \sqrt{41}}{2}.
\]
Since $A$ has integer entries, its characteristic polynomial has rational coefficients. Hence the irrational conjugates $r$ and $s$ must occur with equal multiplicities. 
Since $\dim \mathbf 1^\perp=36$, each must occur with multiplicity $18$.

But then the sum of the nontrivial eigenvalues equals
\[
18(r+s)=18(-3)=-54.
\]
This is impossible, because $G$ is $18$-regular, so $18$ is one eigenvalue of $A$, and the trace identity $\operatorname{tr}(A)=0$ implies that the sum of the remaining $36$ eigenvalues must be $-18$.

This contradiction shows that no such $37$-vertex graph exists. Hence $R(B_8,B_{10})\le 37$. Together with the known lower bound $R(B_8,B_{10})\ge 37$, we obtain $R(B_8,B_{10})=37$, completing the proof of Theorem~\ref{thm:main}.

\section{Lean Formalization}\label{sec:lean}

In order to verify the proof, we formalized the upper-bound argument in Lean~\cite{lean,mathlib}.
The formalization is available on GitHub: \url{https://github.com/BLL4/RamseyB8B10}.

The formalized statement is that every simple graph on $37$ vertices contains a copy of $B_8$, or its complement contains a copy of $B_{10}$. The Lean code for the statement follows.

\begin{center}
\includegraphics[width=11cm]{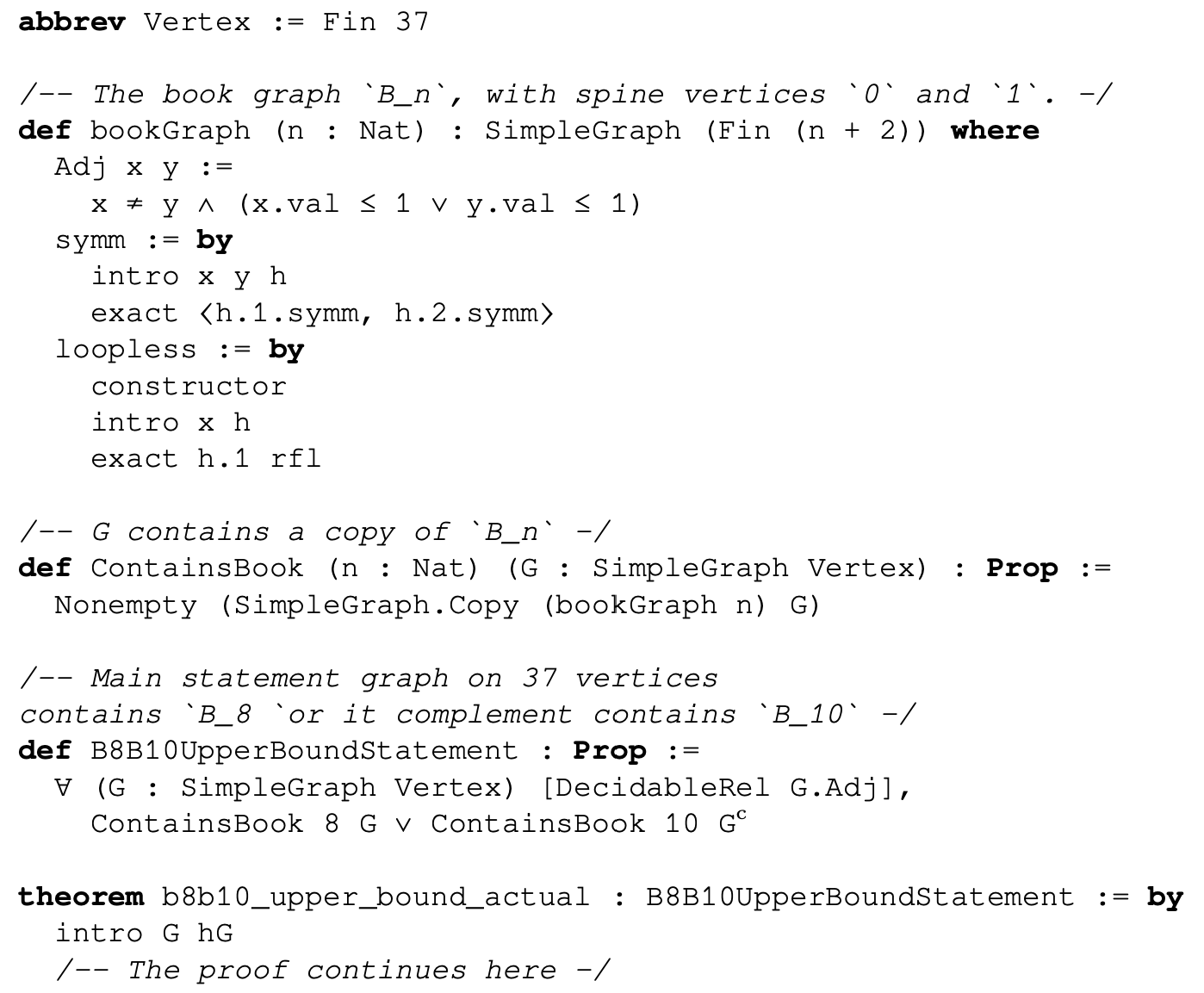}
\end{center}

The Lean formalization was produced with GPT-5.5 Pro through Codex~\cite{codex}.
While the TeX source of this article was provided as an input to the formalization, the proof in Lean is slightly different from the proof presented in the paper in the argument for the non-existence of a strongly regular graph. 

We note that merely asking an AI system to translate the English proof into Lean was not sufficient.
In earlier attempts, AI systems often formalized only the final non-existence of a strongly regular graph with parameters $(37,18,7,10)$, while leaving as an assumption the reduction from the Ramsey hypothesis to those strongly regular parameters.

\section{Note on provenance and software}\label{sec:AI}

The proof of Theorem~\ref{thm:main} was discovered by \emph{AutoMath}~\cite{AutoMath}. 
AutoMath is a publicly available software developed by the first author. 
AutoMath is not a new LLM model.
It is built around command-line-operable language-model agents. 
In the run that led to this note, AutoMath primarily used GPT-5.4 Codex~\cite{codex}. 

\emph{AutoMath} keeps a queue of ``curated'' problems to attempt, a ledger of its work, a registry of attempted problems, and files that contain its work and data for each attempt that it (or its user) may refer back to. Depending on the user's specification, \emph{AutoMath} may work on multiple problems in parallel. In a typical run, it curates a queue of candidate open problems, chooses one or more to attempt, records attempts and data in that problem's respective folder, and then passes successful outputs down to separate verification and audit stages. Lean formalization is an optional part of the \emph{AutoMath} process. When formalization is enabled, solved and verified packets can be put into a separate Lean queue in order to allow \emph{AutoMath} to work on formalization and discovery simultaneously.

For the present theorem, AutoMath was run for approximately 30 minutes (this time includes curation but not Lean formalization, which failed) without parallelization and surfaced the Goodman-identity, regularity, and strongly-regular-graph route; the argument above is the human-written version.

Recent AI-assisted mathematical discovery has taken several forms. 
FunSearch searches over programs using an evaluator and found improved cap-set constructions and online bin-packing heuristics~\cite{FunSearch}. PatternBoost alternates classical local search with transformer-generated seeds to find extremal constructions~\cite{PatternBoost}.
A successor of PatternBoost is Axplorer~\cite{axplorer}.
AlphaEvolve is somewhat an evolution of FunSearch. 
It uses language-model-generated code changes inside an evolutionary evaluation loop and has reported new algorithmic and mathematical constructions~\cite{AlphaEvolve}. In AlphaEvolve, the research team supplies a description of the problem in English and Gemini does the rest.

There is also a line of model-generated proof output for specified mathematical questions. Alexeev, Putterman, Sawhney, Sellke, and Valiant attribute several short proofs in combinatorics, probability, and number theory to an internal OpenAI model~\cite{alexeev2026shortproofscombinatoricsnumber,alexeev2026shortproofscombinatoricsprobability}. Most notably, in 2026, OpenAI announced an internal general-purpose reasoning model's counterexample to Erd\H{o}s's planar unit-distance conjecture~\cite{OpenAIUnitDistance} (see~\cite{UnitDistance} for a shorter, human-verified version). These proofs involve internal model systems and mathematical tasks already selected for investigation.

Conversely, \emph{AutoMath} differs from most AI-assisted proof tools. It is public software rather than an internal model, and it is not ``fed'' problems; rather, it continuously searches for open problems to solve. As such, it can be run without a human in the loop choosing each problem or steering each proof attempt (granted, the final mathematical claim still requires human checking).

\section*{Acknowledgments}

The authors are very grateful to Konstantin Slutsky for his invaluable work on the Lean formalization of the proof.

Konstantin Slutsky pointed out to us that the proof can be simplified by noticing that a strongly regular graph with parameters
$(v,k,\lambda,\mu)=(37,18,7,10)$ is not a conference graph and violates known conditions \cite[p.~222]{Godsil2001-ph}. 
We opted to keep the proof in the current form since we suspect that it may be easier to formalize a self-contained proof in Lean.

\bibliographystyle{plainurl}
\bibliography{references.bib}

\end{document}